# ON THE SIGNAL-TO-INTERFERENCE RATIO OF CDMA SYSTEMS IN WIRELESS COMMUNICATIONS

By Z. D. Bai[1] and Jack W. Silverstein[2]

*Northeast Normal University and National University of Singapore and North Carolina State University*

Let $\{s_{ij}: i,j=1,2,\dots\}$ consist of i.i.d. random variables in $\mathbb{C}$ with $\mathsf{E}s_{11}=0$, $\mathsf{E}|s_{11}|^2=1$. For each positive integer $N$, let $\mathbf{s}_k=\mathbf{s}_k(N)=(s_{1k},s_{2k},\dots,s_{Nk})^T$, $1\le k\le K$, with $K=K(N)$ and $K/N\to c>0$ as $N\to\infty$. Assume for fixed positive integer $L$, for each $N$ and $k\le K$, $\boldsymbol{\alpha}_k=(\alpha_k(1),\dots,\alpha_k(L))^T$ is random, independent of the $s_{ij}$, and the empirical distribution of $(\alpha_1,\dots,\alpha_K)$, with probability one converging weakly to a probability distribution $H$ on $\mathbb{C}^L$. Let $\boldsymbol{\beta}_k=\boldsymbol{\beta}_k(N)=(\alpha_k(1)\mathbf{s}_k^T,\dots,\alpha_k(L)\mathbf{s}_k^T)^T$ and set $C=C(N)=(1/N)\sum_{k=2}^K\boldsymbol{\beta}_k\boldsymbol{\beta}_k^*$. Let $\sigma^2>0$ be arbitrary. Then define $\mathrm{SIR}_1=(1/N)\boldsymbol{\beta}_1^*(C+\sigma^2 I)^{-1}\boldsymbol{\beta}_1$, which represents the best signal-to-interference ratio for user 1 with respect to the other $K-1$ users in a direct-sequence code-division multiple-access system in wireless communications. In this paper it is proven that, with probability 1, $\mathrm{SIR}_1$ tends, as $N\to\infty$, to the limit $\sum_{\ell,\ell'=1}^L\bar\alpha_1(\ell)\alpha_1(\ell')a_{\ell,\ell'}$, where $A=(a_{\ell,\ell'})$ is nonrandom, Hermitian positive definite, and is the unique matrix of such type satisfying $A=(c\,\mathsf{E}\frac{\boldsymbol{\alpha}\boldsymbol{\alpha}^*}{1+\boldsymbol{\alpha}^*A\boldsymbol{\alpha}}+\sigma^2 I_L)^{-1}$, where $\boldsymbol{\alpha}\in\mathbb{C}^L$ has distribution $H$. The result generalizes those previously derived under more restricted assumptions.

**1. Introduction.** This paper examines the mathematical properties of a quantity fundamental in analyzing the performance of a particular scheme used in wireless communications. The scheme, known as *direct-sequence code-division multiple-access* (or DS-CDMA), currently in use, effectively handles many users by taking into account the manner interference interacts when a particular user's information is being decoded. It is achieved by

Received January 2006; revised June 2006.

[1]Supported in part by NSFC Grant 10571020 and NUS Grant R-155-000-056-112.

[2]Supported by the U.S. Army Research Office under Grant W911NF-05-1-0244.

assigning to each user a vector of high dimension, called a *signature sequence*. Suppose there are $K$ users and $L$ receive antennas. Let $N$ be the dimension of the signature sequences, and denote by $\mathbf{s}_k \in \mathbb{C}^N$ the signature sequence assigned to user $k$. At a particular instant of time let $X_k \in \mathbb{R}$ denote the value transmitted by user $k$ having transmit power $T_k \in \mathbb{R}^+$, and let $\gamma_k(\ell)$ denote the fading channel gain from user $k$ to antenna $\ell$. It is assumed that the $X_k$'s are independent standardized random variables. With $\mathbf{W}(\ell) \in \mathbb{C}^N$ denoting noise associated with transmission to antenna $\ell$, entries $W_i(\ell)$ i.i.d. across $i$ and $\ell$, mean zero and $\mathsf{E}|W_i(\ell)|^2 = \sigma^2$, the data recorded at antenna $\ell$ is modeled by

$$\mathbf{Y}(\ell) = \sum_{k=1}^{K} X_k \sqrt{T_k} \gamma_k(\ell) \mathbf{s}_k + \mathbf{W}(\ell).$$

Letting $\mathbf{Y} = [\mathbf{Y}(1)^T, \ldots, \mathbf{Y}(L)^T]^T \in \mathbb{C}^{NL}$, the goal is to capture the transmitted $X_k$ for each user in a linear fashion, that is, by taking the inner product of $\mathbf{Y}$ with an appropriate vector $\mathbf{c}_k \in \mathbb{C}^{NL}$, called the linear receiver for user $k$. For user 1, $\hat{X}_1 = \mathbf{c}_1^* \mathbf{Y}$ is the estimate of transmitted $X_1$.

The output *signal-to-interference ratio*

$$\frac{|\mathbf{c}_1^* \hat{\mathbf{s}}_1|^2}{\sigma^2 \|\mathbf{c}_1\|^2 + \sum_{k=2}^{K} |\mathbf{c}_1^* \hat{\mathbf{s}}_k|^2}$$

associated with user 1 is typically used as a measure for evaluating the performance of the linear receiver. Here

$$\hat{\mathbf{s}}_k = \sqrt{T_k}[\gamma_k(1)\mathbf{s}_k^T, \ldots, \gamma_k(L)\mathbf{s}_k^T]^T.$$

It turns out that the choice of $\mathbf{c}_1$ which minimizes $\mathsf{E}(\hat{X} - X)^2$ (the minimum mean-square error) also maximizes user 1's signal-to-interference ratio, the latter taking the value

$$\mathrm{SIR}_1 = \hat{\mathbf{s}}_1^* \left( \sum_{k=2}^{K} \hat{\mathbf{s}}_k \hat{\mathbf{s}}_k^* + \sigma^2 I \right)^{-1} \hat{\mathbf{s}}_1,$$

where $I$ is the $NL \times NL$ identity matrix. It is this quantity which is the focus of this paper.

In [4] properties of $\mathrm{SIR}_1$ and their dependency on the $\gamma_k(\ell)$'s, $T_k$'s, $\sigma^2$, $L$, $N$ and $K$, when the latter two values are large, are explored by proving limiting results, as $N$ and $K$ approach infinity with their ratio approaching a positive constant, under the assumption that the $s_k$'s are randomly generated (which is usually done in practice). They are independent i.i.d. random vectors containing i.i.d. mean zero entries, independent of the $\gamma_k(\ell)$'s and $T_k$'s. The results allow for analysis of performance in various situations depending on the location of the users with respect to each other and the antennas. Additional assumptions are imposed on the $\gamma_k(\ell)$'s and the $T_k$'s. Throughout [4] it is assumed the $\gamma_k(\ell)$'s are independent and circularly symmetric (i.e., the argument of each $\gamma_k(\ell)$ is uniformly distributed on $[0, 2\pi)$), and the entries of each $\mathbf{s}_k$ are mean zero complex Gaussian with variance $1/N$ [i.e., they are of the form $Z_1 + iZ_2$ with $Z_1$, $Z_2$ i.i.d. $N(0, 1/(2N))$]. The $T_k$'s are allowed only to depend on $|\gamma_j(\ell)|$ for all $k, j, \ell$. Two scenarios depending on the location of the antennas are considered. One scenario places all the antennas near each other, the other allowing them to be located anywhere. Theorem 1 in [4] applies to the former case, the proof of which requires, for each $k$, $\gamma_k(1), \ldots, \gamma_k(L)$ to be identically distributed. Also, $T_k$, as a function of the $|\gamma_j(\ell)|$'s, is assumed to be symmetric with respect to the antennas, in the following sense. For any permutation $\pi$ on $\{1, \ldots, L\}$, we have for each $k$,

$$T_k(\Gamma(\pi(1)), \ldots, \Gamma(\pi(L))) = T_k(\Gamma(1), \ldots, \Gamma(L)),$$

where $\Gamma(\ell) = (\gamma_1(\ell), \ldots, \gamma_k(\ell))$. These two assumptions are lifted in Theorem 3, which would not be realistic when antennas are not placed in one location. However, an additional assumption is made, namely, that there are for each user independent signature sequences going to the $L$ antennas, that is, $\hat{\mathbf{s}}_k$ takes the form

$$\sqrt{T_k}[\gamma_k(1)\mathbf{s}_k^{(1)^T}, \ldots, \gamma_k(L)\mathbf{s}_k^{(L)^T}]^T,$$

with $\mathbf{s}_k^{(1)}, \ldots, \mathbf{s}_k^{(L)}$ i.i.d. As pointed out in [4], this "completely random sequence model is not physically realizable."

The purpose of this paper is to prove limiting results on $\mathrm{SIR}_1$ substantially more general than those found in [4]. The main result is the following:

THEOREM 1.1. *Let $\{s_{ij} : i, j = 1, 2, \ldots\}$ be a doubly infinite array of i.i.d. complex random variables with $\mathsf{E}s_{11} = 0$, $\mathsf{E}|s_{11}|^2 = 1$ (we will from henceforth call standardized). Define for $k = 1, 2, \ldots, K$ $\mathbf{s}_k = \mathbf{s}_k(N) = (s_{1k}, s_{2k}, \ldots, s_{Nk})^T$. We assume $K = K(N)$ and $K/N \to c > 0$ as $N \to \infty$. For each $N$, let $\gamma_k(\ell) = \gamma_k^N(\ell) \in \mathbb{C}$, $T_k = T_k^N \in \mathbb{R}^+$, $k = 1, \ldots, K$, $\ell = 1, \ldots, L$, be random variables, independent of $\mathbf{s}_1, \ldots, \mathbf{s}_K$'s. Let for each $N$ and $k$,*

$$\boldsymbol{\alpha}_k = \boldsymbol{\alpha}_k^N = \sqrt{T_k}(\gamma_k(1), \ldots, \gamma_k(L))^T.$$

*Assume, almost surely, the empirical distribution of $\boldsymbol{\alpha}_1, \ldots, \boldsymbol{\alpha}_K$ weakly converges to a probability distribution $H$ in $\mathbb{C}^L$.*

*Let $\boldsymbol{\beta}_k = \boldsymbol{\beta}_k(N) = \sqrt{T_k}(\gamma_k(1)\mathbf{s}_k^T, \ldots, \gamma_k(L)\mathbf{s}_k^T)^T$ and*

$$C = C(N) = \frac{1}{N}\sum_{k=2}^{K} \boldsymbol{\beta}_k \boldsymbol{\beta}_k^*.$$

*Define*

$$\mathrm{SIR}_1 = \frac{1}{N}\boldsymbol{\beta}_1^*(C+\sigma^2 I)^{-1}\boldsymbol{\beta}_1,$$

*then, with probability one,*

$$\lim_{N\to\infty}\mathrm{SIR}_1 = T_1\sum_{\ell,\ell'=1}^{L}\bar{\gamma}_1(\ell)\gamma_1(\ell')a_{\ell,\ell'},$$

*where the* $L\times L$ *matrix* $A=(a_{\ell,\ell'})$ *is nonrandom, Hermitian positive definite, and is the unique Hermitian positive definite matrix satisfying*

$$A = \left(c\,\mathsf{E}\frac{\boldsymbol{\alpha}\boldsymbol{\alpha}^*}{1+\boldsymbol{\alpha}^*A\boldsymbol{\alpha}}+\sigma^2 I_L\right)^{-1}, \tag{1.1}$$

*where* $\boldsymbol{\alpha}\in\mathbb{C}^L$ *has distribution* $H$ *and* $I_L$ *is the* $L\times L$ *identity matrix.*

Clearly $\mathrm{SIR}_1$ defined in this theorem is the same as the one initially introduced, the only difference in notation being the removal of the scaling by $1/\sqrt{N}$ in the definition of the $\mathbf{s}_k$'s.

Let $\alpha_\ell$ denote the $\ell$th entry of the random vector $\boldsymbol{\alpha}$ having distribution $H$. Under the independence and circularly symmetric assumption on the $\gamma_k(\ell)$'s and the independence of their angles and the $T_k$'s, it follows that for $\ell\neq\ell'$ and positive $a_1,\dots,a_L$,

$$\mathsf{E}\frac{\alpha_\ell\bar{\alpha}_{\ell'}}{1+\sum_{\underline{\ell}}a_{\underline{\ell}}|\alpha_{\underline{\ell}}|^2}=0. \tag{1.2}$$

With just this additional condition we have the following:

Corollary 1.1. *Under the conditions in Theorem 1.1 and (1.2), the limiting* $A=\mathrm{diag}(a_1,\dots,a_L)$*, where the* $a_\ell$*'s are positive satisfying*

$$a_\ell = \frac{1}{c\,\mathsf{E}(|\alpha_\ell|^2/(1+\sum_{\underline{\ell}}a_{\underline{\ell}}|\alpha_{\underline{\ell}}|^2))+\sigma^2}. \tag{1.3}$$

Corollary 1.2. *Suppose the conditions in Theorem 1.1 are met except, for the limiting behavior of the* $\boldsymbol{\alpha}_k$*'s, it is only known that:*

1. *the empirical distribution of*

$$T_k(|\gamma_k(1)|^2,\dots,|\gamma_k(L)|^2)^T,\qquad 2\le k\le K, \tag{1.4}$$

*converges almost surely in distribution to a probability distribution* $G$ *in* $\mathbb{R}^L$*, and*

2. *for $\ell \neq \ell'$ and positive $a_1, \dots, a_L$,*

$$\frac{1}{K-1}\sum_{k=2}^{K}\frac{T_k\gamma_k(\ell)\bar{\gamma}_k(\ell')}{1+\sum_{\underline{\ell}} a_{\underline{\ell}}T_k|\gamma_k(\underline{\ell})|^2} \to 0 \tag{1.5}$$

*almost surely, as $N \to \infty$.*

*Let $(\delta_1, \dots, \delta_L)^T \in \mathbb{R}^L$ denote a random vector having distribution $G$. Then the conclusions of Theorem 1.1 and Corollary 1.1 hold, with each $|\alpha_\ell|^2$ in (1.3) replaced by $\delta_\ell$.*

Theorem 1.1 frees up conditions on the $\gamma_k(\ell)$'s, and overall dependence between them, the $T_k$'s, and the antennas. Moreover, the result allows for more general (and realistic) assumptions on the generation of the $\mathbf{s}_k$'s, permitting their entries, for example, to be just $\pm 1$, which is typically done in practice.

Thus, under the general assumptions in Theorem 1.1, various scenarios can be analyzed and compared. In applications the empirical distribution of the $\boldsymbol{\alpha}_k$'s would typically be used for $H$, the matrix $A$ thereby satisfying

$$A = \left(\frac{K}{N}\frac{1}{K-1}\sum_{k=2}^{K}\frac{\boldsymbol{\alpha}_k\boldsymbol{\alpha}_k^*}{1+\boldsymbol{\alpha}_k^* A\boldsymbol{\alpha}_k} + \sigma^2 I_L\right)^{-1}.$$

Although there appears to be no explicit solution to (1.1), it will be shown that $A$ can be computed numerically by iteration of the right-hand side of (1.1), provided the eigenvalues of the initial choice in the iteration lie in a certain closed interval in $(0, \infty)$.

The conclusion of Corollary 1.1 is the same as that of Theorem 3 in [4], where only the a.s. convergence in distribution of the empirical distribution of (1.4) is assumed. But recall Theorem 3 also assumes for each user different signature sequences for each antenna. The extent of the results in [4] are confined to diagonal $A$'s due to the assumptions imposed [essentially the spherically symmetric assumption on the $\gamma_k(\ell)'s$], clearly a special case of the general conditions assumed in Theorem 1.1.

It is remarked here that the assumption of the $s_{ij}$ coming from a doubly infinite array can be replaced with $s_{ij} = s_{ij}(N)$, $1 \le i \le N$, $1 \le j \le K$, with no dependency assumptions for different $N$, provided $\mathsf{E}|s_{11}|^4 < \infty$. Indeed, it will be seen in the beginning of the proof of Theorem 1.1 that the double array and finite second moment assumption is needed only when the strong law of large numbers is invoked on sums involving $|s_{ij}|^2$, the alternative assumptions yielding the same conclusions with the aid of Lemma 2.10 below.

Theorem 1.1 only provides limiting properties of the signal-to-interference ratio with respect to one user. The last section of this paper will address the issue of uniform convergence of all the $K$ SIR's.

The proofs of these results will be given in Sections 3–7, with basic mathematical results needed in the proofs presented in Section 2.

*Note.* After submitting this paper, the authors came upon a result similar to Theorem 1.1, announced in a conference paper, without proof [3]. In that paper it is claimed that a proof is given in another paper, submitted for publication. In [3] the $s_{ij}$'s need not be identically distributed, nor come from one doubly infinite array of variables, but it is assumed $\mathsf{E}|s_{ij}|^4 = \mathsf{E}|s_{ij}(N)|^4 < N^{2-\gamma}$ for some $\gamma > 1$. Moreover, the limiting distribution $H$ is assumed to have bounded support. The conclusion has convergence in mean square.

**2. Basic tools.** This section contains properties of matrices, a classic fixed point theorem, and some probabilistic results, needed in the proof of the above statements. Throughout, $I$ will denote the $NL \times NL$ identity matrix. For arbitrary dimension $n$, $I_n$ will denote the $n \times n$ identity matrix. For any rectangular matrix $X$, $\operatorname{vec} X$ will denote the column vector consisting of stacking the columns of $X$ on top of each other, first column on top, last on bottom. Spectral norm on matrices and Euclidean norm on vectors will be denoted by $\|\cdot\|$.

LEMMA 2.1. *Let $\sigma^2 > 0$, $B$, $A$ $n \times n$ matrices with $B$ Hermitian nonnegative definite, and $\mathbf{x} \in \mathbb{C}^n$. Then*

$$|\operatorname{tr}((B + \mathbf{x}\mathbf{x}^* + \sigma^2 I)^{-1} - (B + \sigma^2 I)^{-1})A| = \left| \frac{\mathbf{x}^*(B + \sigma^2 I)^{-1} A (B + \sigma^2 I)^{-1}\mathbf{x}}{1 + \mathbf{x}^*(B + \sigma^2 I)^{-1}\mathbf{x}} \right| \le \frac{\|A\|}{\sigma^2}.$$

PROOF. The identity follows from $(D + \mathbf{x}\mathbf{x}^*)^{-1}\mathbf{x} = D^{-1}\mathbf{x}\frac{1}{1+\mathbf{x}^* D^{-1}\mathbf{x}}$, true whenever $n \times n$ $D$ and $D + \mathbf{x}\mathbf{x}^*$ are both invertible. Write $B = \sum \lambda_i \mathbf{e}_i \mathbf{e}_i^*$, its spectral decomposition. Then

$$\begin{aligned} &\left| \frac{\mathbf{x}^*(B + \sigma^2 I)^{-1} A (B + \sigma^2 I)^{-1}\mathbf{x}}{1 + \mathbf{x}^*(B + \sigma^2)^{-1}\mathbf{x}} \right| \\ &\quad \le \frac{\|A\| \, \|(B + \sigma^2 I)^{-1}\mathbf{x}\|^2}{1 + \mathbf{x}^*(B + \sigma^2 I)^{-1}\mathbf{x}} \\ &\quad = \|A\| \frac{\sum (1/(\lambda_i + \sigma^2)^2)|\mathbf{e}_i^*\mathbf{x}|^2}{1 + \sum (1/(\lambda_i + \sigma^2))|\mathbf{e}_i^*\mathbf{x}|^2} \\ &\quad \le \frac{\|A\|}{\sigma^2}. \end{aligned}$$

□

The next lemma is easily verifiable.

LEMMA 2.2. *For any matrix $A$ $N\times K$ and $\sigma^2>0$,*

$$(AA^*+\sigma^2 I_N)^{-1}=\sigma^{-2}(I_N-A(A^*A+\sigma^2 I_K)^{-1}A^*).$$

The following lemma is a direct consequence of the previous identity.

LEMMA 2.3. *Suppose $A_1,\dots,A_L$ are $N\times K$, and $\sigma^2>0$. Define the $\ell,\ell'$ block of the $NL\times NL$ matrix $\mathcal{A}$ by $\mathcal{A}_{\ell,\ell'}=A_\ell A^*_{\ell'}$ and, splitting $(\mathcal{A}+\sigma^2 I)^{-1}$ into $L^2$ $N\times N$ matrices, let $(\mathcal{A}+\sigma^2 I)^{-1}_{\ell,\ell'}$ denote its $\ell,\ell'$ block. Then*

$$(\mathcal{A}+\sigma^2 I)^{-1}_{\ell,\ell'}=\sigma^{-2}\left(\delta_{\ell,\ell'}I_N-A_\ell\left(\sum_{\underline{\ell}}A^*_{\underline{\ell}}A_{\underline{\ell}}+\sigma^2 I_K\right)^{-1}A^*_{\ell'}\right).$$

LEMMA 2.4. *Given $A_1,\dots,A_L$ are $N\times K$ and $z_1,\dots,z_\ell\in\mathbb{C}$ with $\sum_\ell |z_\ell|^2=1$,*

$$\left(\sum_\ell A_\ell z_\ell\right)\left(\sum_\ell A^*_\ell \hat{z}_\ell\right)\preceq\sum_\ell A_\ell A^*_\ell,$$

*where "$\preceq$" represents the partial ordering on Hermitian nonnegative definite matrices.*

PROOF. For any $\mathbf{b}\in\mathbb{C}^N$, we have by two applications of Cauchy–Schwarz,

$$\begin{aligned}\mathbf{b}^*\left(\sum_\ell A_\ell z_\ell\right)\left(\sum_\ell A^*_\ell \hat{z}_\ell\right)\mathbf{b}&=\sum_{\ell,\ell'}\mathbf{b}^*A_\ell A^*_{\ell'}\mathbf{b}z_\ell\hat{z}_\ell\\&\le\sum_{\ell,\ell'}(\mathbf{b}^*A_\ell A^*_\ell\mathbf{b})^{1/2}(\mathbf{b}^*A_{\ell'}A^*_{\ell'}\mathbf{b})^{1/2}|z_\ell|\,|z_{\ell'}|\\&=\left(\sum_\ell(\mathbf{b}^*A_\ell A^*_\ell\mathbf{b})^{1/2}|z_\ell|\right)^2\\&\le\sum_\ell|z_\ell|^2\sum_\ell\mathbf{b}^*A_\ell A^*_\ell\mathbf{b}=\mathbf{b}^*\left(\sum_\ell A_\ell A^*_\ell\right)\mathbf{b}.\end{aligned}$$

This proves the result. □

LEMMA 2.5. *For $A_1,\dots,A_L$, $\mathcal{A}$, $\sigma^2$ in Lemma 2.3, the $L\times L$ matrix $(\operatorname{tr}(\mathcal{A}+\sigma^2 I)^{-1}_{\ell,\ell'})$ is positive definite with smallest eigenvalue bounded below by*

$$\operatorname{tr}\left(\sum_\ell A_\ell A^*_\ell+\sigma^2 I_N\right)^{-1}.$$

PROOF. For $z_1,\dots,z_L\in\mathbb{C}$, $\sum_\ell|z_\ell|^2=1$, we have by Lemmas 2.2–2.4 the smallest eigenvalue of $(\operatorname{tr}(A+\sigma^2I)^{-1}_{\ell,\ell'})$ is bounded below by

$$\begin{aligned}
&\sum_{\ell,\ell'}\operatorname{tr}(A+\sigma^2I)^{-1}_{\ell,\ell'}\bar z_\ell\bar z_{\ell'}\\
&\quad=\sigma^{-2}\operatorname{tr}\left(\sum_{\ell,\ell'}\delta_{\ell,\ell'}\bar z_\ell z_{\ell'}I_N-\sum_{\ell,\ell'}A_\ell\left(\sum_\ell A_\ell^*A_\ell+\sigma^2I_K\right)^{-1}A_{\ell'}\bar z_\ell z_{\ell'}\right)\\
&\quad=\sigma^{-2}\operatorname{tr}\left(I_N-\left(\sum_\ell A_\ell z_\ell\right)\left(\sum_\ell A_\ell^*A_\ell+\sigma^2I_K\right)^{-1}\left(\sum_\ell A_\ell^*\bar z_\ell\right)\right)\\
&\quad\ge\sigma^{-2}\operatorname{tr}\left(I_N-\left(\sum_\ell A_\ell z_\ell\right)\right.\\
&\qquad\qquad\left.\times\left(\left(\sum_\ell A_\ell^*\bar z_\ell\right)\left(\sum_\ell A_\ell z_\ell\right)+\sigma^2I_K\right)^{-1}\left(\sum_\ell A_\ell^*\bar z_\ell\right)\right)\\
&\quad=\operatorname{tr}\left(\left(\sum_\ell A_\ell z_\ell\right)\left(\sum_\ell A_\ell^*\bar z_\ell\right)+\sigma^2I_N\right)^{-1}\ge\operatorname{tr}\left(\sum_\ell A_\ell A_\ell^*+\sigma^2I_N\right)^{-1}.\qquad\square
\end{aligned}$$

For $A=(a_{ij})$ $m\times n$ and $B$ $p\times q$, the Kronecker product of $A$ and $B$, denoted by $A\otimes B$, is the $mp\times nq$ matrix, expressed in blocks of $p\times q$ matrices, the $i,j$ block being $a_{ij}B$. We will need the following, which is Lemma 4.2.10 of [5].

LEMMA 2.6. *For $A$ $m\times n$, $B$ $p\times q$, $C$ $n\times k$ and $D$ $q\times r$, we have*

$$(A\otimes B)(C\otimes D)=(AC)\otimes(BD).$$

The following is needed to prove Corollary 1.1.

LEMMA 2.7 (Schauder fixed point theorem [7]). *If $\mathcal{A}$ is a convex, compact subset of a Banach space $\mathcal{X}$ and $g:\mathcal{A}\to\mathcal{A}$ is continuous, then $g$ has a fixed point in $\mathcal{A}$.*

The next result is one on the eigenvalues of the expected value of the Kronecker product of two random matrices.

LEMMA 2.8. *Let $A=(a_{ij})=(\mathbf{a}_1,\dots,\mathbf{a}_n)$ $(m\times n)$ and $B$ $(h\times g)$ be two random matrices, the entries having bounded second moments. Then*

$$\|\mathsf{E}A\otimes B\|\le\min\left(\sqrt{\|\mathsf{E}AA^*\|\,\|\mathsf{E}B^*B\|},\sqrt{\|\mathsf{E}A^*A\|\,\|\mathsf{E}BB^*\|}\right).$$

PROOF. For any $h\times m$ $X=(\mathbf{x}_1,\dots,\mathbf{x}_m)$ and $g\times n$ matrix $Y=(\mathbf{y}_1,\dots,\mathbf{y}_n)$ with $\operatorname{tr} XX^*=\operatorname{tr} YY^*=1$, we have, using Cauchy–Schwarz,

$$\begin{aligned}
&\|(\mathrm{vec}X)^*[\mathsf{E}A\otimes B]\mathrm{vec}Y\|^2\\
&\quad=\left|\sum_{i=1}^m\sum_{j=1}^n \mathsf{E}a_{ij}\mathbf{x}_i^*B\mathbf{y}_j\right|^2=\left|\sum_{j=1}^n \mathsf{E}\mathbf{a}_j^TX^*B\mathbf{y}_j\right|^2\\
&\quad\le\left(\sum_{j=1}^n \mathsf{E}[\|\overline{X}\mathbf{a}_j\|\,\|B\mathbf{y}_j\|]\right)^2\le\left(\sum_{j=1}^n \mathsf{E}^{1/2}\|\overline{X}\mathbf{a}_j\|^2\mathsf{E}^{1/2}\|B\mathbf{y}_j\|^2\right)^2\\
&\quad\le\sum_{j=1}^n \mathsf{E}\|\overline{X}\mathbf{a}_j\|^2\sum_{j=1}^n \mathsf{E}\|B\mathbf{y}_j\|^2=\operatorname{tr}(\mathsf{E}AA^*\overline{X}^*\overline{X})\operatorname{tr}(\mathsf{E}B^*BYY^*)\\
&\quad\le\|\mathsf{E}AA^*\|\,\|\mathsf{E}B^*B\|.
\end{aligned}$$

Notice that

$$\sum_{j=1}^n \mathsf{E}\mathbf{a}_j^TX^*B\mathbf{y}_j=\operatorname{tr}\mathsf{E}A^TX^*BY=\operatorname{tr}\mathsf{E}B^T\overline{X}AY^T,$$

so we also have

$$\|(\mathrm{vec}X)^*[\mathsf{E}A\otimes B]\mathrm{vec}Y\|^2\le\|\mathsf{E}A^*A\|\,\|\mathsf{E}BB^*\|.$$

The truth of the lemma follows. □

The next result, which is Lemma 2.7 in [2], constitutes the main contribution of randomness to Theorem 1.1.

LEMMA 2.9. *For $\mathbf{X}=(X_1,\dots,X_n)^T$ i.i.d. standardized entries, $C$ $n\times n$, we have for any $p\ge 2$,*

$$\mathsf{E}|\mathbf{X}^*C\mathbf{X}-\operatorname{tr}C|^p\le K_p((\mathsf{E}|X_1|^4\operatorname{tr}CC^*)^{p/2}+\mathsf{E}|X_1|^{2p}\operatorname{tr}(CC^*)^{p/2}),$$

*where the constant $K_p$ does not depend on $n$, $C$, nor on the distribution of $X_1$.*

The last two results provide conditions guaranteeing the strong law of large numbers.

LEMMA 2.10. ([6]). *For $X_1,X_2,\dots$ i.i.d., let $S_n=X_1+\cdots+X_n$. For $t\ge 1$, the joint conditions $\mathsf{E}|X_1|^t<\infty$ and $\mathsf{E}X_1=b$ are equivalent to the condition*

$$\sum_{n=1}^\infty n^{t-2}\mathsf{P}\left(\left|\frac{S_n}{n}-b\right|\ge\varepsilon\right)<\infty$$

*for every $\varepsilon>0$.*

LEMMA 2.11 (Lemma 2 of [1]). *Let $\{X_{ij},\, i,j=1,2,\ldots\}$ be a double array of i.i.d. random variables and let $\alpha > \frac{2}{3}$, $\beta \ge 0$ and $M > 0$ be constants. Then as $n \to \infty$,*

$$\max_{j \le Mn^\beta} \left| n^{-\alpha} \sum_{i=1}^{n} (x_{ij} - c) \right| \to 0 \qquad a.s.$$

*if and only if the following hold:*

$$\mathsf{E}|X_{11}|^{(1+\beta)/\alpha} < \infty$$

*and*

$$c = \begin{cases} \mathsf{E}X_{11}, & \text{if } \alpha \le 1, \\ \text{any number}, & \text{if } \alpha > 1. \end{cases}$$

**3. Proof of Theorem 1.1.** For the remainder of this paper we write $\sqrt{T_k}\gamma_k(\ell)$ as $\alpha_k(\ell)$. Write

$$\mathrm{SIR}_1 = \frac{1}{N} \sum_{\ell,\ell'} \bar{\alpha}_1(\ell)\alpha_1(\ell')\mathbf{s}_1^*(C + \sigma^2 I)^{-1}_{\ell,\ell'}\mathbf{s}_1,$$

where $(C+\sigma^2 I)^{-1}_{\ell,\ell'}$, $N \times N$, is the $\ell,\ell'$ block of the $NL \times NL$ matrix $(C + \sigma^2 I)^{-1}$. Some of the $NL \times NL$ matrices below will also be viewed in block form: $L^2$ $N \times N$ matrices.

We begin by truncating and centralizing the entries of $\mathbf{s}_1$. For each $N$, define

$$\tilde{s}_{n1}(N) = s_{n1} I_{(|s_{n1}| \le (1/3)\log N)} - \mathsf{E}s_{n1} I_{(|s_{n1}| \le (1/3)\log N)},$$

$$\tilde{\mathbf{s}}_1 = \tilde{\mathbf{s}}_1(N) = (\tilde{s}_{11}(N), \tilde{s}_{21}(N), \ldots, \tilde{s}_{N1}(N))^T,$$

$$\tilde{\boldsymbol{\beta}}_1 = \tilde{\boldsymbol{\beta}}_1(N) = (\alpha_1(1)\tilde{\mathbf{s}}_1^T, \ldots, \alpha_1(L)\tilde{\mathbf{s}}_1)^T,$$

and $\widetilde{\mathrm{SIR}}_1 = \frac{1}{N}\tilde{\boldsymbol{\beta}}_1^*(C+\sigma^2 I)^{-1}\tilde{\boldsymbol{\beta}}_1 = \frac{1}{N}\sum_{\ell,\ell'} \bar{\alpha}_1(\ell)\alpha_1(\ell')\tilde{\mathbf{s}}_1^*(C+\sigma^2 I)^{-1}_{\ell,\ell'}\tilde{\mathbf{s}}_1$.

We have

$$\begin{aligned} |\mathrm{SIR}_1 - \widetilde{\mathrm{SIR}}_1| &\le \sigma^{-2}N^{-1}(\|\boldsymbol{\beta}_1 - \tilde{\boldsymbol{\beta}}_1\|^2 + 2\|\boldsymbol{\beta}_1\|\,\|\boldsymbol{\beta}_1 - \tilde{\boldsymbol{\beta}}_1\|) \\ &= \sigma^{-2}\sum_{\ell} |\alpha_1(\ell)|^2 N^{-1}(\|\mathbf{s}_1 - \tilde{\mathbf{s}}_1\|^2 + 2\|\mathbf{s}_1\|\|\mathbf{s}_1 - \tilde{\mathbf{s}}_1\|). \end{aligned}$$

We have by the strong law of large numbers (SLLN) $N^{-1}\|\mathbf{s}_1\|^2 \to 1$ a.s. as $N \to \infty$. Also, since

$$\begin{aligned} N^{-1}\|\mathbf{s}_1 - \tilde{\mathbf{s}}_1\|^2 &= N^{-1}\sum_{n=1}^{N} |s_{n1} I_{(|s_{n1}| > (1/3)\log N)} - \mathsf{E}s_{n1} I_{(|s_{n1}| \le (1/3)\log N)}|^2 \\ &\le N^{-1} 2 \sum_{n=1}^{N} |s_{n1}|^2 I_{(|s_{n1}| > (1/3)\log N)} + 2|\mathsf{E}s_{11} I_{(|s_{11}| \le (1/3)\log N)}|^2, \end{aligned}$$

we have, for any positive $M$ by SLLN,

$$\limsup_N N^{-1}\|\mathbf{s}_1-\tilde{\mathbf{s}}_1\|^2 \le \lim_{N\to\infty} N^{-1}\sum_{n=1}^N 2|s_{n1}I_{(|s_{n1}|>M)}|^2$$

$$= 2\mathsf{E}|s_{11}|^2 I_{(|s_{11}|>M)} \qquad \text{a.s.}$$

Since this last expression can be made arbitrarily small by choosing $M$ sufficiently large, we have

$$|\mathrm{SIR}_1 - \widetilde{\mathrm{SIR}}_1| \to 0 \qquad \text{a.s.}$$

Let $\hat{\mathbf{s}}_1 = \hat{\mathbf{s}}_1(N) = \tilde{\mathbf{s}}_1(N)/(\mathsf{E}|\tilde{s}_{11}(N)|^2)^{1/2}$. It is clear that $\mathsf{E}|\tilde{s}_{11}(N)|^2 \to 1$ as $N\to\infty$. Therefore,

$$\left|\widetilde{\mathrm{SIR}}_1 - \frac{1}{N}\sum_{\ell,\ell'}\bar{\alpha}_1(\ell)\alpha_1(\ell')\hat{\mathbf{s}}_1^*(C+\sigma^2 I)^{-1}_{\ell,\ell'}\hat{\mathbf{s}}_1\right| \to 0$$

as $N\to\infty$. Since the entries of $\hat{\mathbf{s}}_1$ are i.i.d. standardized and, for all large $N$, bounded by $\log N$, we have by Lemma 2.9, for any $N\times N$ $A$,

$$\mathsf{E}|N^{-1}\hat{\mathbf{s}}_1^* A\hat{\mathbf{s}}_1 - N^{-1}\operatorname{tr} A|^4 \le K_4\|A\|^4(\log N)^8 N^{-2},$$

which is summable for $A$ bounded in norm. Since we have $\|(C+\sigma^2 I)^{-1}_{\ell,\ell'}\| \le \sigma^{-2}$, we conclude that, with probability one,

$$\left|\mathrm{SIR}_1 - N^{-1}\sum_{\ell,\ell'}\bar{\alpha}_1(\ell)\alpha_1(\ell')\operatorname{tr}(C+\sigma^2 I)^{-1}_{\ell,\ell'}\right| \to 0.$$

Similar to what is argued above, we have by SLLN

$$\lim_{N\to\infty}\frac{1}{N^2}\sum_{n=1}^N\sum_{k=2}^K |s_{nk}|^2 I_{(|s_{nk}|>\log N)} = 0 \qquad \text{a.s.}$$

It is straightforward to verify the existence of a nonrandom sequence $\{a_N\}$, of positive numbers increasing to infinity, satisfying

$$\lim_{N\to\infty}\frac{a_N^2}{N^2}\sum_{n=1}^N\sum_{k=2}^K |s_{nk}|^2 I_{(|s_{nk}|>\log N)} = 0 \tag{3.1}$$

almost surely, and

$$a_N \mathsf{E} s_{11} I_{(|s_{11}|\le \log N)} \to 0. \tag{3.2}$$

We may assume $a_N \le \log N$.

Define $\mathcal{K}_1 = \mathcal{K}(N) = \{k\in\{2,\dots,K\} : \max_\ell |\alpha_k(\ell)|^2 < a_N\}$ and $\mathcal{K}_2 = \{2,\dots,K\} - \mathcal{K}_1$. Since the empirical distribution of $\boldsymbol{\alpha}_k$ converges weakly a.s. to a probability distribution in $\mathbb{C}^L$, we must have, with probability one,

$$\#\mathcal{K}_2 \equiv K_N = o(N).$$

Write $\mathcal{K}_2 = \{k_1, \dots, k_{K_N}\}$, $\hat{C}_0 = C$, and $\hat{C}_j = \hat{C}_{j-1} - \boldsymbol{\beta}_{k_j}\boldsymbol{\beta}^*_{k_j}$, $j = 1, \dots, K_N$. Let $I_{\ell',\ell}$ denote the $NL \times NL$ matrix consisting of the $N \times N$ identity matrix in the $\ell', \ell$ block, zeros elsewhere. Then, using Lemma 2.1 for any $\ell, \ell'$,

$$
\begin{aligned}
&|N^{-1}\operatorname{tr}(C+\sigma^2)^{-1}_{\ell,\ell'} - N^{-1}\operatorname{tr}(\hat{C}_{K_N} + \sigma^2 I)^{-1}_{\ell,\ell'}| \\
&\quad = \left|\sum_{j=1}^{K_N} \operatorname{tr}(\hat{C}_{j-1} + \sigma^2 I)^{-1} I_{\ell'\ell} - N^{-1}\operatorname{tr}(\hat{C}_j + \sigma^2 I)^{-1} I_{\ell'\ell}\right| \\
&\quad = \left|\frac{1}{N}\sum_{j=1}^{K_N} \frac{N^{-1}\boldsymbol{\beta}^*_{k_j}(\hat{C}_j + \sigma^2 I)^{-1} I_{\ell',\ell}(\hat{C}_j + \sigma^2 I)^{-1}\boldsymbol{\beta}_{k_j}}{1 + N^{-1}\boldsymbol{\beta}^*_{k_j}(\hat{C}_j + \sigma^2 I)^{-1}\boldsymbol{\beta}_{k_j}}\right| \\
&\quad \le \frac{K_N}{N\sigma^2} \to 0 \qquad \text{a.s.}
\end{aligned}
$$

Therefore, we may assume each $\boldsymbol{\alpha}_{k_j} = 0$, and that $\max_{2\le k\le K,\ell} |\boldsymbol{\alpha}_k(\ell)|^2 \le a_N$.

We truncate and centralize the entries of each $\mathbf{s}_k$, $2 \le k \le K$, in the same manner as performed on $\mathbf{s}_1$, and call it $\tilde{\mathbf{s}}_k$. The corresponding $\boldsymbol{\beta}_k$ and $C$ are denoted by $\tilde{\boldsymbol{\beta}}_k$ and $\tilde{C}$. For any $\ell, \ell'$, we have

$$
\begin{aligned}
&|N^{-1}(\operatorname{tr}(C + \sigma^2 I)^{-1}_{\ell,\ell'} - \operatorname{tr}(\tilde{C} + \sigma^2 I)^{-1}_{\ell,\ell'})| \\
&\quad = \left|\frac{1}{N^2}\sum_{k=2}^{K} \operatorname{tr}((C+\sigma^2 I)^{-1}(\boldsymbol{\beta}_k\boldsymbol{\beta}^*_k - \tilde{\boldsymbol{\beta}}_k\tilde{\boldsymbol{\beta}}^*_k)(\tilde{C} + \sigma^2 I)^{-1} I_{\ell',\ell})\right| \\
&\quad = \left|\frac{1}{N^2}\sum_{k=2}^{K} \boldsymbol{\beta}^*_k(\tilde{C} + \sigma^2 I)^{-1} I_{\ell',\ell}(C + \sigma^2 I)^{-1}\boldsymbol{\beta}_k \right. \\
&\qquad\qquad \left. - \tilde{\boldsymbol{\beta}}^*_k(\tilde{C} + \sigma^2 I)^{-1} I_{\ell',\ell}(C + \sigma^2 I)^{-1}\tilde{\boldsymbol{\beta}}_k\right| \\
&\quad \le \frac{1}{N^2\sigma^4}\sum_{k=2}^{K} \|\boldsymbol{\beta}_k - \tilde{\boldsymbol{\beta}}_k\|^2 + 2\|\boldsymbol{\beta}_k\|\|\boldsymbol{\beta}_k - \tilde{\boldsymbol{\beta}}_k\| \\
&\quad = \frac{1}{N^2\sigma^4}\sum_{k=2}^{K}\sum_{\ell} |\alpha_k(\ell)|^2(\|\mathbf{s}_k - \tilde{\mathbf{s}}_k\|^2 + 2\|\mathbf{s}_k\|\|\mathbf{s}_k - \tilde{\mathbf{s}}_k\|) \\
&\quad \le \frac{a_N L}{N^2\sigma^4}\left[\sum_{k=2}^{K} \|\mathbf{s}_k - \tilde{\mathbf{s}}_k\|^2 + 2\left(\sum_{k=2}^{K}\|\mathbf{s}_k\|^2\right)^{1/2}\left(\sum_{k=2}^{K}\|\mathbf{s}_k - \tilde{\mathbf{s}}_k\|^2\right)^{1/2}\right].
\end{aligned}
$$

By SLLN, we have

$$
\frac{1}{NK}\sum_{k=2}^{K}\|\mathbf{s}_k\|^2 = \frac{1}{NK}\sum_{n=1}^{N}\sum_{k=2}^{K}|s_{nk}|^2 \to 1 \qquad \text{a.s.}
$$

and from (3.1) and (3.2),

$$\frac{a_N^2}{N^2}\sum_{k=2}^K \|\mathbf{s}_k - \tilde{\mathbf{s}}_k\|^2$$
$$\le \frac{2a_N^2}{NK}\sum_{n=1}^N\sum_{k=2}^K |s_{nk}|^2 I_{(|s_{nk}|>(1/3)\log N)}$$
$$+ 2a_N^2|\mathsf{E}s_{11}I_{(|s_{11}|\le(1/3)\log N)}|^2 \to 0 \qquad \text{a.s.}$$

Therefore,

$$N^{-1}(\mathrm{tr}(C+\sigma^2 I)^{-1}_{\ell,\ell'} - \mathrm{tr}(\tilde{C}+\sigma^2 I)^{-1}_{\ell,\ell'}) \to 0 \qquad \text{a.s.}$$

It is easy to verify

$$N^{-1}(\mathrm{tr}((\mathsf{E}|\tilde{s}_{11}|^2)^{-1}\tilde{C}+\sigma^2 I)^{-1}_{\ell,\ell'} - \mathrm{tr}(\tilde{C}+\sigma^2 I)^{-1}_{\ell,\ell'}) \to 0.$$

Therefore, returning to the original notation, we may replace the doubly infinite array and assume for each $N$, $s_{nk}=s_{nk}(N)$, $1\le n\le N$, $2\le k\le K$, i.i.d. standardized random variables bounded by $\log N$. The quantities $\mathbf{s}_k=\mathbf{s}_k(N)$, $\boldsymbol{\beta}_k=\boldsymbol{\beta}_k(N)$, and $C=C(N)$ are defined accordingly.

Define $C_{(k)}=C-(1/N)\boldsymbol{\beta}_k\boldsymbol{\beta}_k^*$. Select $\varepsilon\in(0,1/10)$. Applying Lemma 2.9 to each of $\mathbf{s}_k$, $(C_{(k)}+\sigma^2)^{-1}_{\ell,\ell'}$, $2\le k\le K$, $1\le \ell,\ell'\le L$, with $p=5$, along with standard arguments using Chebyshev's and Boole's inequalities, together with Lemma 2.1, we have

$$\max_{\substack{k\in\{2,\dots,K\}\\ \ell,\ell'\in\{1,\dots,L\}}}\left|\frac{N^\varepsilon}{N}(\mathbf{s}_k^*(C_{(k)}+\sigma^2 I)^{-1}_{\ell,\ell'}\mathbf{s}_k - \mathrm{tr}(C+\sigma^2 I)^{-1}_{\ell,\ell'})\right| \to 0 \qquad \text{a.s.} \tag{3.3}$$

Define the $L\times L$ matrix $\underline{B}=(\underline{b}_{\ell,\ell'})$ with

$$\underline{b}_{\ell,\ell'} = \frac{1}{N}\sum_{k=2}^K \frac{\bar{\alpha}_k(\ell')\alpha_k(\ell)}{1+(1/N)\boldsymbol{\beta}_k^*(C_{(k)}+\sigma^2 I)^{-1}\boldsymbol{\beta}_k},$$

and define the $NL\times NL$ matrix $B$ in terms of the Kronecker product: $B=\underline{B}\otimes I_N$. We have $(B+\sigma^2 I)^{-1} = (\underline{B}+\sigma^2 I_L)^{-1}\otimes I_N$. Denote the $\ell,\ell'$ entry of $(\underline{B}+\sigma^2 I_L)^{-1}$ by $\hat{b}_{\ell,\ell'}$.

We write

$$C+\sigma^2 I - (B+\sigma^2 I) = \frac{1}{N}\sum_{k=2}^K \boldsymbol{\beta}_k\boldsymbol{\beta}_k^* - B.$$

Taking inverses on each side, we have

$$(B+\sigma^2 I)^{-1} - (C+\sigma^2 I)^{-1}$$

$$= \frac{1}{N}\sum_{k=2}^{K}(B+\sigma^2 I)^{-1}\boldsymbol{\beta}_k\boldsymbol{\beta}_k^*(C+\sigma^2 I)^{-1} - (B+\sigma^2 I)^{-1}B(C+\sigma^2 I)^{-1}$$

$$= \frac{1}{N}\sum_{k=2}^{K}\frac{(B+\sigma^2 I)^{-1}\boldsymbol{\beta}_k\boldsymbol{\beta}_k^*(C_{(k)}+\sigma^2 I)^{-1}}{1+(1/N)\boldsymbol{\beta}_k^*(C_{(k)}+\sigma^2 I)^{-1}\boldsymbol{\beta}_k} - (B+\sigma^2 I)^{-1}B(C+\sigma^2 I)^{-1}.$$

Multiplying on the right by $I_{\ell',\ell}$, taking traces and dividing by $N$, we get

$$\begin{aligned}
&N^{-1}\operatorname{tr}(B+\sigma^2 I)^{-1}_{\ell,\ell'} - N^{-1}\operatorname{tr}(C+\sigma^2 I)^{-1}_{\ell,\ell'}\\
&\quad= \frac{1}{N}\sum_{k=2}^{K}\frac{(1/N)\boldsymbol{\beta}_k^*(C_{(k)}+\sigma^2 I)^{-1}I_{\ell',\ell}(B+\sigma^2 I)^{-1}\boldsymbol{\beta}_k}{1+(1/N)\boldsymbol{\beta}_k^*(C_{(k)}+\sigma^2 I)^{-1}\boldsymbol{\beta}_k}\\
&\qquad - N^{-1}\operatorname{tr}B(C+\sigma^2 I)^{-1}I_{\ell',\ell}(B+\sigma^2 I)^{-1}\\
&\quad= \frac{1}{N}\sum_{k=2}^{K}\frac{\sum_{\underline{\ell},\underline{\ell}'}\bar{\alpha}_k(\underline{\ell})\alpha_k(\underline{\ell}')(1/N)\mathbf{s}_k^*[(C_{(k)}+\sigma^2 I)^{-1}I_{\ell',\ell}(B+\sigma^2 I)^{-1}]_{\underline{\ell},\underline{\ell}'}\mathbf{s}_k}{1+(1/N)\boldsymbol{\beta}_k^*(C_{(k)}+\sigma^2 I)^{-1}\boldsymbol{\beta}_k}\\
&\qquad - \frac{1}{N}\operatorname{tr}\sum_{\underline{\ell},\underline{\ell}'}B_{\underline{\ell}',\underline{\ell}}[(C+\sigma^2 I)^{-1}I_{\ell',\ell}(B+\sigma^2 I)^{-1}]_{\underline{\ell},\underline{\ell}'}\\
&\quad= \sum_{\underline{\ell},\underline{\ell}'}\frac{1}{N}\Bigg[\sum_{k=2}^{K}\frac{1}{N}\mathbf{s}_k^*[(C_{(k)}+\sigma^2 I)^{-1}I_{\ell',\ell}(B+\sigma^2 I)^{-1}]_{\underline{\ell},\underline{\ell}'}\mathbf{s}_k\\
&\qquad\qquad\times\frac{\bar{\alpha}_k(\underline{\ell})\alpha_k(\underline{\ell}')}{1+(1/N)\boldsymbol{\beta}_k^*(C_{(k)}+\sigma^2 I)^{-1}\boldsymbol{\beta}_k}\\
&\qquad\qquad-\operatorname{tr}B_{\underline{\ell}',\underline{\ell}}[(C+\sigma^2 I)^{-1}I_{\ell',\ell}(B+\sigma^2 I)^{-1}]_{\underline{\ell},\underline{\ell}'}\Bigg]\\
&\quad= \sum_{\underline{\ell},\underline{\ell}'}\frac{1}{N}\Bigg[\sum_{k=2}^{K}\frac{1}{N}\mathbf{s}_k^*(C_{(k)}+\sigma^2 I)^{-1}_{\underline{\ell},\ell'}(B+\sigma^2 I)^{-1}_{\ell,\underline{\ell}'}\mathbf{s}_k\\
&\qquad\qquad\times\frac{\bar{\alpha}_k(\underline{\ell})\alpha_k(\underline{\ell}')}{1+(1/N)\boldsymbol{\beta}_k^*(C_{(k)}+\sigma^2 I)^{-1}\boldsymbol{\beta}_k}\\
&\qquad\qquad-\operatorname{tr}B_{\underline{\ell}',\underline{\ell}}(C+\sigma^2 I)^{-1}_{\underline{\ell},\ell'}(B+\sigma^2 I)^{-1}_{\ell,\underline{\ell}'}\Bigg]\\
&\quad= \sum_{\underline{\ell},\underline{\ell}'}\frac{1}{N}\sum_{k=2}^{K}\frac{\hat{b}_{\ell,\underline{\ell}'}\bar{\alpha}_k(\underline{\ell})\alpha_k(\underline{\ell}')}{1+(1/N)\boldsymbol{\beta}_k^*(C_{(k)}+\sigma^2 I)^{-1}\boldsymbol{\beta}_k}N^{-1}\\
&\qquad\qquad\times(\mathbf{s}_k^*(C_{(k)}+\sigma^2 I)^{-1}_{\underline{\ell},\ell'}\mathbf{s}_k - \operatorname{tr}(C+\sigma^2 I)^{-1}_{\underline{\ell},\ell'}).
\end{aligned}$$

Using (3.3), the fact that the $\hat{b}_{\ell,\underline{\ell'}}$'s are bounded by $\sigma^{-2}$, and noticing that

$$N^{-1}\operatorname{tr}(B+\sigma^2 I)^{-1}_{\ell,\ell'} = \hat{b}_{\ell,\ell'},$$

we immediately get

$$|\hat{b}_{\ell,\ell'} - N^{-1}\operatorname{tr}(C+\sigma^2 I)^{-1}_{\ell,\ell'}| \to 0 \qquad \text{a.s.} \tag{3.4}$$

Notice $C_{\ell,\ell'}$, the $\ell,\ell'$ block of $C$, can be written as

$$((1/\sqrt{N})S\mathcal{A}(\ell))((1/\sqrt{N})S\mathcal{A}(\ell'))^*,$$

where $S = (\mathbf{s}_2,\ldots,\mathbf{s}_K)$ and $\mathcal{A}(\ell) = \operatorname{diag}(\alpha_2(\ell),\ldots,\alpha_k(\ell))$. Therefore, from Lemma 2.5, the $L\times L$ matrix $(\operatorname{tr}(C+\sigma^2 I)^{-1}_{\ell,\ell'})$ is positive definite.

Using this and (3.3), we have

$$\begin{aligned}&\left|\underline{b}_{\ell,\ell'} - \frac{1}{N}\sum_{k=2}^{K}\frac{\bar{\alpha}_k(\ell')\alpha_k(\ell)}{1+\sum_{\underline{\ell},\underline{\ell'}}\bar{\alpha}_k(\underline{\ell})\alpha_k(\underline{\ell'})N^{-1}\operatorname{tr}(C+\sigma^2 I)^{-1}_{\underline{\ell},\underline{\ell'}}}\right| \\ &\quad\le \sum_{\underline{\ell},\underline{\ell'}}\frac{a_N^2}{N}\sum_{k=2}^{N}N^{-1}|\mathbf{s}_k^*(C_{(k)}+\sigma^2 I)^{-1}_{\underline{\ell},\underline{\ell'}}\mathbf{s}_k \\ &\qquad\qquad - \operatorname{tr}(C+\sigma^2 I)^{-1}_{\underline{\ell},\underline{\ell'}}| \to 0 \qquad \text{a.s.}\end{aligned} \tag{3.5}$$

From Lemma 2.5 the smallest eigenvalue of $A_N \equiv (N^{-1}\operatorname{tr}(C+\sigma^2 I)^{-1}_{\ell,\ell'})$ is bounded below by

$$\frac{1}{N}\operatorname{tr}\left(\frac{1}{N}S\sum_{\ell}\mathcal{A}(\ell)\mathcal{A}(\ell)^* S^* + \sigma^2 I_N\right)^{-1}. \tag{3.6}$$

This quantity is the *Stieltjes transform* of the empirical distribution of the eigenvalues of

$$\frac{1}{N}S\sum_{\ell}\mathcal{A}(\ell)\mathcal{A}(\ell)^* S^*$$

evaluated at $-\sigma^2$. We have, with probability one, the empirical distribution of the diagonal entries of

$$\sum_{\ell}\mathcal{A}(\ell)\mathcal{A}(\ell)^*$$

converging weakly to a nonrandom probability distribution. Therefore, from [8], we see, with probability one, the empirical distribution of the eigenvalues of (3.6) converges weakly to a nonrandom probability distribution, and consequently, (3.6) converges a.s. to a nonrandom positive number, say, $m$. Therefore,

$$\liminf_{N}\lambda_{\min}A_N \ge m \qquad \text{a.s.} \tag{3.7}$$

Consider a realization in which (3.4),(3.5) and (3.7) hold and the empirical distribution of $(\boldsymbol{\alpha}_2,\dots,\boldsymbol{\alpha}_K)$ converges weakly, where $\boldsymbol{\alpha}\in\mathbb{C}^L$ denotes the random vector having distribution $H$. Let $\{N_i\}$ be a subsequence in which each $N^{-1}\operatorname{tr}(C+\sigma^2 I)^{-1}_{\ell,\ell'}$ converges, say, to $a_{\ell,\ell'}$ for $\ell,\ell'\in\{1,\dots,L\}$. Let $A=(a_{\ell,\ell'})$ and $\delta=\inf_{N_i}\lambda_{\min}(A_{N_i})>0$. We have for $\mathbf{z}\in\mathbb{C}^L$,

$$\left|\frac{\mathbf{z}\mathbf{z}^*}{1+\mathbf{z}^*A\mathbf{z}}\right|\le\frac{\|\mathbf{z}\|^2}{1+\lambda_{\min}(A)\|\mathbf{z}\|^2}\le\frac{1}{\delta}.$$

Therefore, by the dominated convergence theorem, along $\{N_i\}$,

$$\frac{1}{K-1}\sum_{k=2}^{K}\frac{\boldsymbol{\alpha}_k\boldsymbol{\alpha}_k^*}{1+\boldsymbol{\alpha}_k^*A\boldsymbol{\alpha}_k}\to\mathsf{E}\frac{\boldsymbol{\alpha}\boldsymbol{\alpha}^*}{1+\boldsymbol{\alpha}^*A\boldsymbol{\alpha}},$$

and since

$$\left|\frac{\mathbf{z}\mathbf{z}^*}{1+\mathbf{z}^*A_N\mathbf{z}}-\frac{\mathbf{z}\mathbf{z}^*}{1+\mathbf{z}^*A\mathbf{z}}\right|\le\frac{1}{\delta}\|A_N-A\|,$$

we have by (3.5), along $\{N_i\}$,

$$\underline{B}\to c\,\mathsf{E}\frac{\boldsymbol{\alpha}\boldsymbol{\alpha}^*}{1+\boldsymbol{\alpha}^*A\boldsymbol{\alpha}}.$$

So $A$ satisfies (1.1). The next section shows that only one Hermitian positive definite $A$ will satisfy this equation. With this fact we have, with probability one, $A_N$ converges to a nonrandom Hermitian positive definite $L\times L$ matrix $A$ satisfying (1.1).

**4. Proof of uniqueness.** Suppose $A$ and $\widetilde{A}$ are two different $L\times L$ Hermitian positive definite matrices satisfying (1.1). Then

$$A-\widetilde{A}=c\,\mathsf{E}\frac{A\boldsymbol{\alpha}\boldsymbol{\alpha}^*\widetilde{A}\boldsymbol{\alpha}^*(A-\widetilde{A})\boldsymbol{\alpha}}{(1+\boldsymbol{\alpha}^*A\boldsymbol{\alpha})(1+\boldsymbol{\alpha}^*\widetilde{A}\boldsymbol{\alpha})}.$$

Multiplying $A^{-1/2}$ on the left and $\widetilde{A}^{-1/2}$ on the right, we obtain

$$\begin{aligned}A^{1/2}\widetilde{A}^{-1/2}-A^{-1/2}\widetilde{A}^{1/2}&=c\mathsf{E}\frac{A^{1/2}\boldsymbol{\alpha}\boldsymbol{\alpha}^*\widetilde{A}^{1/2}\boldsymbol{\alpha}^*(A-\widetilde{A})\boldsymbol{\alpha}}{(1+\boldsymbol{\alpha}^*A\boldsymbol{\alpha})(1+\boldsymbol{\alpha}^*\widetilde{A}\boldsymbol{\alpha})}\\&=c\mathsf{E}\frac{\boldsymbol{\eta}\widetilde{\boldsymbol{\eta}}^*\boldsymbol{\eta}^*(A^{1/2}\widetilde{A}^{-1/2}-A^{-1/2}\widetilde{A}^{1/2})\widetilde{\boldsymbol{\eta}}}{(1+\boldsymbol{\alpha}^*A\boldsymbol{\alpha})(1+\boldsymbol{\alpha}^*\widetilde{A}\boldsymbol{\alpha})},\end{aligned}$$

where $\boldsymbol{\eta}=A^{1/2}\boldsymbol{\alpha}$ and $\widetilde{\boldsymbol{\eta}}=\widetilde{A}^{1/2}\boldsymbol{\alpha}$. Write $\boldsymbol{\mu}=\operatorname{vec}(A^{1/2}\widetilde{A}^{-1/2}-A^{-1/2}\widetilde{A}^{1/2})$. With the aid of the Kronecker product, we can write the above equation as

$$\boldsymbol{\mu}=c\mathsf{E}\frac{(\overline{\widetilde{\boldsymbol{\eta}}}\otimes\boldsymbol{\eta})(\widetilde{\boldsymbol{\eta}}^T\otimes\boldsymbol{\eta}^*)}{(1+\boldsymbol{\alpha}^*A\boldsymbol{\alpha})(1+\boldsymbol{\alpha}^*\widetilde{A}\boldsymbol{\alpha})}\boldsymbol{\mu}.\tag{4.1}$$

Using Lemma 2.6, we have

$$c\mathsf{E}\frac{(\overline{\widetilde{\boldsymbol{\eta}}}\otimes\boldsymbol{\eta})(\widetilde{\boldsymbol{\eta}}^T\otimes\boldsymbol{\eta}^*)}{(1+\boldsymbol{\alpha}^*A\boldsymbol{\alpha})(1+\boldsymbol{\alpha}^*\widetilde{A}\boldsymbol{\alpha})}=c\mathsf{E}\left[\frac{\overline{\widetilde{\boldsymbol{\eta}}\widetilde{\boldsymbol{\eta}}^*}}{1+\boldsymbol{\alpha}^*\widetilde{A}\boldsymbol{\alpha}}\otimes\frac{\boldsymbol{\eta}\boldsymbol{\eta}^*}{1+\boldsymbol{\alpha}^*\widetilde{A}\boldsymbol{\alpha}}\right]$$

and, since $\boldsymbol{\mu}\neq 0$, this matrix has an eigenvalue equal to 1. By Lemma 2.8, its largest squared eigenvalue cannot be greater than

$$\left\|c\mathsf{E}\left(\frac{\widetilde{\boldsymbol{\eta}}\widetilde{\boldsymbol{\eta}}^*}{1+\boldsymbol{\alpha}^*\widetilde{A}\boldsymbol{\alpha}}\right)^2\right\|\left\|c\mathsf{E}\left(\frac{\boldsymbol{\eta}\boldsymbol{\eta}^*}{1+\boldsymbol{\alpha}^*A\boldsymbol{\alpha}}\right)^2\right\|.$$

We have

$$c\mathsf{E}\left(\frac{\boldsymbol{\eta}\boldsymbol{\eta}^*}{1+\boldsymbol{\alpha}^*A\boldsymbol{\alpha}}\right)^2=c\mathsf{E}\frac{A^{1/2}\boldsymbol{\alpha}\boldsymbol{\alpha}^*A^{1/2}\boldsymbol{\alpha}^*A\boldsymbol{\alpha}}{(1+\boldsymbol{\alpha}^*A\boldsymbol{\alpha})^2},$$

and since

$$\frac{A^{1/2}\boldsymbol{\alpha}\boldsymbol{\alpha}^*A^{1/2}}{1+\boldsymbol{\alpha}^*A\boldsymbol{\alpha}}-\frac{A^{1/2}\boldsymbol{\alpha}\boldsymbol{\alpha}^*A^{1/2}\boldsymbol{\alpha}^*A\boldsymbol{\alpha}}{(1+\boldsymbol{\alpha}^*A\boldsymbol{\alpha})^2}$$

is nonnegative definite, we have

$$\begin{aligned}c\mathsf{E}\left(\frac{\boldsymbol{\eta}\boldsymbol{\eta}^*}{1+\boldsymbol{\alpha}^*A\boldsymbol{\alpha}}\right)^2&\preceq c\mathsf{E}\frac{A^{1/2}\boldsymbol{\alpha}\boldsymbol{\alpha}^*A^{1/2}}{1+\boldsymbol{\alpha}^*A\boldsymbol{\alpha}}\\&=A^{1/2}(A^{-1}-\sigma^2I_L)A^{1/2}=I_L-\sigma^2A,\end{aligned}$$

the eigenvalues of which must all be less than one. The same result applies for the other matrix involving $\widetilde{A}$. Therefore, the matrix in (4.1) cannot have an eigenvalue equal to one, a contradiction. So we conclude that there is only one Hermitian positive definite solution to (1.1).

## 5. Convergence of iterations.

Let $f(A)$ denote the right-hand side of (1.1), considered as a mapping of the set of Hermitian positive definite matrices, which we will denote by $\mathcal{H}$. Clearly $f$ maps $\mathcal{H}$ into itself with largest eigenvalue not larger than $\sigma^{-2}$. We proceed in finding a positive $b<\sigma^{-2}$ such that $f$ maps

$$\mathcal{H}[b,\sigma^{-2}]\equiv\{A\in\mathcal{H}\colon \text{ all eigenvalues of } A \text{ lie in } [b,\sigma^{-2}]\}$$

into itself. Notice from the dominated convergence theorem

$$g(a)\equiv c\mathsf{E}\frac{a\|\boldsymbol{\alpha}\|^2}{1+a\|\boldsymbol{\alpha}\|^2}+a\sigma^2$$

is continuous and nondecreasing for $a\in[0,\sigma^{-2}]$ with $g(0)=0$ and $g(\sigma^{-2})\geq 1$. Therefore, there exists $\hat{a}\in(0,\sigma^{-2}]$ for which $g(\hat{a})=1$. We claim a suitable $b$ is $\hat{a}/(c+1)$. Indeed, suppose the eigenvalues of $A\in\mathcal{H}$ are contained in $[\hat{a}/(c+1), \sigma^{-2}]$. If $a \equiv \lambda_{\min}(A) \geq \hat{a}$, then using the fact that $\|\mathsf{E}B\| \leq \mathsf{E}\|B\|$ for any random matrix $B$,

$$\begin{aligned}\lambda_{\max}(f^{-1}(A)) &= \sigma^2 + c\lambda_{\max}\left(\mathsf{E}\frac{\boldsymbol{\alpha}\boldsymbol{\alpha}^*}{1+\boldsymbol{\alpha}^* A\boldsymbol{\alpha}}\right)\\ &\leq \sigma^2 + c\mathsf{E}\frac{\|\boldsymbol{\alpha}\|^2}{1+a\|\boldsymbol{\alpha}\|^2} \leq \frac{1+c}{a} \leq \frac{1+c}{\hat{a}},\end{aligned}$$

whereas if $a \in [\hat{a}/(c+1), \hat{a})$,

$$\lambda_{\max} f^{-1}(A) \leq \frac{g(a)}{a} \leq \frac{1}{a} \leq \frac{c+1}{\hat{a}}.$$

The claim is proven.

Let $A_0 \in \mathcal{H}[b, \sigma^{-2}]$ and define recursively $A_{n+1} = f(A_n)$. We have $A_n \in \mathcal{H}[b, \sigma^{-2}]$ for all $n$, and for $n \geq 1$, as in Section 4,

$$A_{n+1}^{1/2} A_n^{-1/2} - A_{n+1}^{-1/2} A_n^{1/2} = c\mathsf{E}\frac{A_{n+1}^{1/2}\boldsymbol{\alpha}\boldsymbol{\alpha}^* A_n^{1/2}\boldsymbol{\alpha}^*(A_n - A_{n-1})\boldsymbol{\alpha}}{(1+\boldsymbol{\alpha}^* A_n\boldsymbol{\alpha})(1+\boldsymbol{\alpha}^* A_{n-1}\boldsymbol{\alpha})}.$$

Letting

$$H_n = A_n^{1/2} A_{n-1}^{-1/2} - A_n^{-1/2} A_{n-1}^{1/2},$$

we have

$$\begin{aligned}H_{n+1} &= c\mathsf{E}\frac{A_{n+1}^{1/2}\boldsymbol{\alpha}\boldsymbol{\alpha}^* A_n^{1/2}\boldsymbol{\alpha}^* A_n^{1/2} H_n A_{n-1}^{1/2}\boldsymbol{\alpha}}{(1+\boldsymbol{\alpha}^* A_n\boldsymbol{\alpha})(1+\boldsymbol{\alpha}^* A_{n-1}\boldsymbol{\alpha})}\\ &= c\mathsf{E}\frac{A_{n+1}^{1/2}\boldsymbol{\alpha}\boldsymbol{\alpha}^* A_n^{1/2}\boldsymbol{\alpha}^* A_{n-1}^{1/2} H_n^* A_n^{1/2}\boldsymbol{\alpha}}{(1+\boldsymbol{\alpha}^* A_n\boldsymbol{\alpha})(1+\boldsymbol{\alpha}^* A_{n-1}\boldsymbol{\alpha})},\end{aligned}$$

or in vector form,

$$\begin{aligned}\mathrm{vec} H_{n+1} &= \left[c\mathsf{E}\frac{(\bar{A}_n^{1/2}\hat{\boldsymbol{\alpha}} \otimes A_{n+1}^{1/2}\boldsymbol{\alpha})(\boldsymbol{\alpha}^T \bar{A}_n^{1/2} \otimes \boldsymbol{\alpha}^* A_{n-1}^{1/2})}{(1+\boldsymbol{\alpha}^* A_n\boldsymbol{\alpha})(1+\boldsymbol{\alpha}^* A_{n-1}\boldsymbol{\alpha})}\right]\mathrm{vec} H_n^*\\ &= c\mathsf{E}\left[\frac{(\bar{A}_n^{1/2}\bar{\boldsymbol{\alpha}}\boldsymbol{\alpha}^T \bar{A}_n^{1/2}) \otimes (A_{n+1}^{1/2}\boldsymbol{\alpha}\boldsymbol{\alpha}^* A_{n-1}^{1/2})}{(1+\boldsymbol{\alpha}^* A_{n-1})(1+\boldsymbol{\alpha}^* A_n\boldsymbol{\alpha})}\right]\mathrm{vec} H_n^*,\end{aligned}$$

using Lemma 2.6. Arguing the same way as in the previous section, we have by Lemma 2.8,

$$\begin{aligned}&\left\|c\mathsf{E}\left[\frac{(\bar{A}_n^{1/2}\bar{\boldsymbol{\alpha}}\boldsymbol{\alpha}^T \bar{A}_n^{1/2}) \otimes (A_{n+1}^{1/2}\boldsymbol{\alpha}\boldsymbol{\alpha}^* A_{n-1}^{1/2})}{(1+\boldsymbol{\alpha}^* A_{n-1})(1+\boldsymbol{\alpha}^* A_n\boldsymbol{\alpha})}\right]\right\|^2\\ &\quad \leq \left\|c\mathsf{E}\frac{A_n^{1/2}\boldsymbol{\alpha}\boldsymbol{\alpha}^* A_n\boldsymbol{\alpha}\boldsymbol{\alpha}^* A_n^{1/2}}{(1+\boldsymbol{\alpha}^* A_{n-1}\boldsymbol{\alpha})(1+\boldsymbol{\alpha}^* A_n\boldsymbol{\alpha})}\right\|\left\|c\mathsf{E}\frac{A_{n+1}^{1/2}\boldsymbol{\alpha}\boldsymbol{\alpha}^* A_{n-1}\boldsymbol{\alpha}\boldsymbol{\alpha}^* A_{n+1}^{1/2}}{(1+\boldsymbol{\alpha}^* A_{n-1}\boldsymbol{\alpha})(1+\boldsymbol{\alpha}^* A_n\boldsymbol{\alpha})}\right\|\end{aligned}$$

$$\le \left\| c\mathsf{E}\frac{A_n^{1/2}\boldsymbol{\alpha}\boldsymbol{\alpha}^* A_n^{1/2}}{1+\boldsymbol{\alpha}^* A_{n-1}\boldsymbol{\alpha}} \right\| \left\| c\mathsf{E}\frac{A_{n+1}^{1/2}\boldsymbol{\alpha}\boldsymbol{\alpha}^* A_{n+1}^{1/2}}{1+\boldsymbol{\alpha}^* A_n\boldsymbol{\alpha}} \right\|$$
$$= \|A_n^{1/2}(A_n^{-1}-\sigma^2 I_L)A_n^{1/2}\| \|A_{n+1}^{1/2}(A_{n+1}^{-1}-\sigma^2 I_L)A_{n+1}^{1/2}\|$$
$$= \|I_L - \sigma^2 A_n\| \, \|I_L - \sigma^2 A_{n+1}\|$$
$$\le (1-\sigma^2 b)^2.$$

Notice $\rho \equiv 1-\sigma^2 b \in (0,1)$. For $n \ge 2$, we therefore get

$$\|H_n\| \le \|\mathrm{vec} H_n\| \le \rho^{n-1}\|\mathrm{vec} H_1\|,$$

and so

$$\|A_n - A_{n-1}\| = \|A_n^{1/2} H_n A_{n-1}^{1/2}\| \le \sigma^{-2}\rho^{n-1}\|\mathrm{vec} H_1\|,$$

which implies for $m \ge n \ge 2$,

$$\|A_m - A_n\| \le \frac{\|\mathrm{vec} H_1\|}{\sigma^2(1-\rho)}\rho^n.$$

Therefore, $\{A_n\}$ is a Cauchy sequence, and hence, convergent to a matrix $A \in \mathcal{H}[b,\sigma^2]$. From continuity of $f$, $A$ satisfies (1.1).

**6. Proofs of corollaries.** For the first corollary we see that, under assumption (1.2), $f$ maps diagonal matrices consisting of positive diagonal elements into diagonal matrices. Due to the uniqueness of solutions to (1.1), the proof amounts to showing the existence of positive $a_1,\dots,a_L$ satisfying (1.1). This is achieved by invoking Lemma 2.7. We simply take $\mathcal{X} = \mathbb{R}^L$, $g$ the right-hand side of (1.3) and $\mathcal{A} = [b,\sigma^{-2}]^L$. The first statement in Corollary 1.1 follows.

For the second corollary, we follow along the argument toward the end of Section 3. We see immediately that (3.7) holds. Consider a realization in which (3.4), (3.5) and (3.7) hold, the empirical distribution of $(|\alpha_k(1)|^2,\dots, |\alpha_k(L)|^2)$, $2 \le k \le K$, converges weakly to $G$, and (1.5) is true for all positive rational $a_1,\dots,a_L$. Then, for this realization, a simple continuity argument reveals (1.5) true for all positive $a_1,\dots,a_L$. Moreover, the empirical distribution of $\boldsymbol{\alpha}_k$, $2 \le k \le K$, is tight. The subsequence $\{N_i\}$ considered can therefore also be one in which the empirical distribution of $\boldsymbol{\alpha}_k$, $2 \le k \le K$, converges weakly to, say, $H$. The rest of the argument at the end of Section 3 leads to only one solution $A$ on $\{N_i\}$ satisfying (1.1). But, by the dominated convergence theorem, (1.2) holds for all positive $a_1,\dots,a_L$. Thus, from Corollary 1.1, $A$ is diagonal satisfying (1.3), which depends only on $G$. Thus, Corollary 1.2 follows.

**7. Question of uniformity.** Let

$$C_k = \frac{1}{N}\left(\sum_{j=1}^{K} \boldsymbol{\beta}_j \boldsymbol{\beta}_j^* - \boldsymbol{\beta}_k \boldsymbol{\beta}_k^*\right).$$

Then

$$\mathrm{SIR}_k \equiv \frac{1}{N}\boldsymbol{\beta}_k^*(C_k+\sigma^2 I)^{-1}\boldsymbol{\beta}_k = \frac{1}{N}\sum_{\ell,\ell'} \bar{\alpha}_k(\ell)\alpha_k(\ell')\mathbf{s}_k^*(C_k+\sigma^2 I)^{-1}_{\ell,\ell'}\mathbf{s}_k$$

represents user $k$'s best signal-to-interference ratio. We are interested in knowing what conditions are needed to insure

$$\max_{k\le K}\left|\frac{1}{N}\sum_{\ell,\ell'} \bar{\alpha}_k(\ell)\alpha_k(\ell')(\mathbf{s}_k^*(C_k+\sigma^2 I)^{-1}_{\ell,\ell'}\mathbf{s}_k - a_{\ell,\ell'})\right| \to \qquad \text{a.s.}$$

Clearly nothing can be concluded without assuming bounds or some growth rate on the $\alpha_k(\ell)$'s along with knowledge of the rate of convergence of the $(1/N)\mathbf{s}_k^*(C_k+\sigma^2 I)^{-1}_{\ell,\ell'}\mathbf{s}_k$'s. The latter is tied closely with its limiting distributional behavior, which will be investigated in later work. For now we will confine the analysis to providing conditions to ensure for any $\ell,\ell'$,

$$\max_{k\le K}|N^{-1}\mathbf{s}_k^*(C_k+\sigma^2 I)^{-1}_{\ell,\ell'}\mathbf{s}_k - a_{\ell,\ell'}| \to 0 \qquad \text{a.s.} \tag{7.1}$$

as $N\to\infty$.

We have the following:

THEOREM 7.1. *If, in addition to the conditions in Theorem 1.1,* $\mathsf{E}|s_{11}|^4 < \infty$*, or if the doubly infinite array assumption is dropped,* $\mathsf{E}|s_{11}|^6 < \infty$*, then (7.1) is true.*

PROOF. For each $k \le K$, let $\tilde{\mathbf{s}}_k$ denote the vector obtained after truncating and centralizing $\mathbf{s}_k$, the same way as $\mathbf{s}_1$. Each $C_k$ remains unchanged. We have

$$\begin{aligned}&|N^{-1}\mathbf{s}_k^*(C_k+\sigma^2 I)^{-1}_{\ell,\ell'}\mathbf{s}_k - N^{-1}\tilde{\mathbf{s}}_k^*(C_k+\sigma^2 I)^{-1}_{\ell,\ell'}\tilde{\mathbf{s}}_k|\\ &\quad\le \sigma^{-2}N^{-1}(\|\mathbf{s}_k-\tilde{\mathbf{s}}_k\|^2 + 2\|\mathbf{s}_k\|\|\mathbf{s}_k-\tilde{\mathbf{s}}_k\|).\end{aligned}$$

By Lemma 2.11, ($X_{11} = |s_{11}|^2$, $\alpha=\beta=1$) under the double array assumption, or, for nondouble array, Lemma 2.10 ($X_1 = |s_{11}|^2$, $t=3$) together with Boole's inequality, we follow the steps in the beginning of Section 3 and find, almost surely,

$$\limsup_N \max_{k\le K}|N^{-1}\|\mathbf{s}_k\|^2 - 1| = 0$$

and

$$\limsup_N \max_{k\le K} N^{-1}\|\mathbf{s}_k-\tilde{\mathbf{s}}_k\|^2=0.$$

Letting $\hat{\mathbf{s}}_k=\tilde{\mathbf{s}}_k/(\mathsf{E}|\tilde{s}_{11}|^2)^{1/2}$, it follows that, almost surely,

$$\max_{k\le K} N^{-1}|\mathbf{s}_k^*(C_k+\sigma^2 I)^{-1}_{\ell,\ell'}\mathbf{s}_k-\hat{\mathbf{s}}_k^*(C_k+\sigma^2 I)^{-1}_{\ell,\ell'}\hat{\mathbf{s}}_k|\to 0$$

as $N\to\infty$. Applying Lemma 2.9 for any $p>4$, we have then

$$\max_{k\le K} N^{-1}|\mathbf{s}_k^*(C_k+\sigma^2 I)^{-1}_{\ell,\ell'}\mathbf{s}_k-\operatorname{tr}(C_k+\sigma^2 I)^{-1}_{\ell,\ell'}|\to 0 \qquad \text{a.s.}$$

as $N\to\infty$. For any $k,k'$, with two applications of Lemma 2.1, we find (with $A=I_{\ell',\ell}$)

$$N^{-1}|\operatorname{tr}(C_k+\sigma^2 I)^{-1}_{\ell,\ell'}-\operatorname{tr}(C_{k'}+\sigma^2 I)^{-1}_{\ell,\ell'}|\le 2\sigma^{-2}N^{-1}.$$

Thus, the remainder of the proof of Theorem 1.1 proceeds exactly as in Section 3, and we get (7.1). $\square$

**Acknowledgments.** The authors wish to thank Brian Hughes, Ajith Kamath and Antonia Tulino for their help in explaining the engineering aspects of the quantity studied. They also thank the referees for their helpful suggestions.

## REFERENCES

[1] Bai, Z. D. and Yin, Y. Q. (1993). Limit of the smallest eigenvalue of a large-dimensional sample covariance matrix. *Ann. Probab.* **21** 1275–1294. MR1235416

[2] Bai, Z. D. and Silverstein, J. W. (1998). No eigenvalues outside the support of the limiting spectral distribution of large-dimensional random matrices. *Ann. Probab.* **26** 316–345. MR1617051

[3] Cottatellucci, L. and Müller, R. (2004). A generalized resourse pooling result for correlated antennas with applications to asychronous CDMA. *International Symposium on Information Theory and Its Applications* (*ISITA 2004*) *Parma, Italy, October 10–13, 2004.*

[4] Hanly, S. V. and Tse, D. N. C. (2001). Resource pooling and effective bandwidths in CDMA networks with multiuser receivers and spatial diversity. *IEEE Trans. Inform. Theory* **47** 1328–1351. MR1830083

[5] Horn, R. A. and Johnson, C. R. (1991). *Topics in Matrix Analysis*. Cambridge Univ. Press. MR1091716

[6] Katz, M. (1963). The probability in the tail of a distribution. *Ann. Math. Statist.* **34** 312–318. MR0144369

[7] Nirenberg, L. (1961). *Functional Analysis*. Lectures given in 1960–61, Courant Institute of Mathematical Sciences. Notes by Lesly Sibner.

[8] Silverstein, J. W. and Bai, Z. D. (1995). On the empirical distribution of eigenvalues of a class of large dimensional random matrices. *J. Multivariate Anal.* **54** 175–192. MR1345534

KLASMOE
School of Mathematics
and Statistics
Northeast Normal University
Changchun, Jilin 130024
PR China
E-mail: baizd@nenu.edu.cn

Department of Mathematics
North Carolina State University
Box 8205
Raleigh, North Carolina 27695-8205
USA
E-mail: jack@math.ncsu.edu